
\documentstyle{amsppt}
\baselineskip18pt
\magnification=\magstep1
\pagewidth{30pc}
\pageheight{45pc}

\hyphenation{co-deter-min-ant co-deter-min-ants pa-ra-met-rised
pre-print pro-pa-gat-ing pro-pa-gate
fel-low-ship Cox-et-er dis-trib-ut-ive}
\def\leaderfill{\leaders\hbox to 1em{\hss.\hss}\hfill}

\

\def\be{{\beta}}

\def\bc{{\bold c}}

\def\bi{{\bold i}}
\def\bj{{\bold j}}
\def\bk{{\bold k}}

\def\br{{\bold r}}

\def\bu{{\bold u}}
\def\bv{{\bold v}}

\def\b0{\text{\bf 0}}

\def\boxit#1{\vbox{\hrule\hbox{\vrule \kern3pt
\vbox{\kern3pt\hbox{#1}\kern3pt}\kern3pt\vrule}\hrule}}
\def\rabbit{\vbox{\hbox{\kern0pt
\vbox{\kern0pt{\hbox{---}}\kern3.5pt}}}}

\def\tableau#1{
        \hbox {
                \hskip -10pt plus0pt minus0pt
                \raise\baselineskip\hbox{
                \offinterlineskip
                \hbox{#1}}
                \hskip0.25em
        }
}

\def\tabCol#1{
\hbox{\vtop{\hrule
\halign{\strut\vrule\hskip0.5em##\hskip0.5em\hfill\vrule\cr\lower0pt
\hbox\bgroup$#1$\egroup \cr}
\hrule
} } \hskip -10.5pt plus0pt minus0pt}

\def\CR{
        $\egroup\cr
        \noalign{\hrule}
        \lower0pt\hbox\bgroup$
}



\def\blank#1#2{
\hbox to #1{\hfill \vbox to #2{\vfill}}
}


\def\strut{\vrule height10pt depth5pt width0pt}

\def\secy{1}
\def\secza{2}
\def\seczc{3}
\def\seczd{4}
\def\secze{5}
\def\seca{6}
\def\secb{7}
\def\al{\alpha}
\def\ga{\gamma}

\topmatter
\title On the maximally clustered elements of Coxeter groups
\endtitle

\author R.M. Green \endauthor
\affil Department of Mathematics \\ University of Colorado \\
Campus Box 395 \\ Boulder, CO  80309-0395 \\ USA \\ {\it  E-mail:}
rmg\@euclid.colorado.edu \\
\newline
\endaffil

\abstract 
We continue the study of the maximally clustered elements
for simply laced Coxeter groups which were recently introduced by Losonczy.
Such elements include as a
special case the freely braided elements of Losonczy and the author,
which in turn constitute a superset of the $iji$-avoiding
elements of Fan.  Our main result is to classify the MC-finite Coxeter
groups, namely those Coxeter groups having
finitely many maximally clustered elements.  Remarkably, any simply laced 
Coxeter group having finitely many $iji$-avoiding elements also turns out to 
be MC-finite.
\endabstract


 \endkeywords

\endtopmatter

\centerline{\bf Preliminary version, draft 2}

\subhead \secy. Introduction \endsubhead

Let $W$ be a simply laced Coxeter group with set $S$ of Coxeter
generators.  Recently, Losonczy \cite{{\bf 7}} introduced
the notion of maximally clustered elements for $W$, and studied
applications of maximally clustered elements to Schubert varieties.

Maximally clustered elements are defined in terms of certain triples of
root vectors.   We explain this briefly as follows. Let $w \in W$.
Every reduced expression for $w$ determines a sequence whose terms are
the positive roots sent negative by $w$.  If such a ``root sequence" has
a consecutive subsequence of the form $\al, \al + \be, \be$, then the set
containing these three vectors is called a ``contractible triple" of $w$.
Some previously studied classes of group elements can be characterized in
terms of contractible triples. For example, the $iji$-avoiding elements of
Fan \cite{{\bf 3}} are precisely those elements of $W$ having no contractible
triples, and the freely braided elements of the author and
Losonczy \cite{{\bf 4}, {\bf 5}} are the elements of $W$ with pairwise disjoint 
contractible triples.

A maximally clustered element $w \in W$ is one with the property that if $T$
and $T'$ are contractible triples of $w$ and $T \cap T' \neq
\emptyset$, then the highest roots of $T$ and $T'$ agree.  Maximally
clustered elements exist in abundance: for example, in the Coxeter group of
type $A_3$, $21$ of the $24$ elements are maximally clustered.
Maximally clustered elements have convenient reduced expressions called
``contracted reduced expressions'', and every reduced expression is short braid
equivalent (commutation equivalent) to a contracted one.  A main result of
\cite{{\bf 7}} is a criterion for a Schubert variety indexed by a maximally 
clustered element to be smooth, and this is much simpler than the corresponding
situation for an arbitrary group element.

The main result of this paper is a classification of simply laced Coxeter
groups $W$ having finitely many maximally clustered elements; we
call these MC-finite Coxeter groups for short.  It is not 
hard to show that for simply laced Coxeter groups, the maximally clustered 
elements are a superset
of the freely braided elements, and that in turn, the freely braided elements
are a superset of the fully commutative elements of Stembridge \cite{{\bf 9}},
or equivalently the $iji$-avoiding elements of Fan \cite{{\bf 3}}.  In order
for $W$ to be MC-finite, it is therefore
necessary for $W$ to have finitely many fully commutative elements.
Our main result is that this condition is also sufficient to
ensure that $W$ is MC-finite.
(It is possible for $W$ to be both infinite and MC-finite.)

Our main result generalizes one of the two main results in \cite{{\bf 5}}, which
proves that a simply laced Coxeter group has finitely many fully commutative
elements if and only if it has finitely many freely braided elements.
Remarkably, the proof presented here is simpler than the proof of the less
general result in \cite{{\bf 5}}, although the argument of this paper 
relies on the classification of Coxeter groups having finitely many fully
commutative elements, whereas the argument of \cite{{\bf 5}} is conceptual.

The paper is organized as follows.  Sections \secza--\seczd\  recall the
necessary background from the theory of Coxeter groups.  The main result
is stated in Section \secze.  Section \seca\  develops the theory of certain
operators, $\pi_i$, on reduced expressions, and these are used to prove the
main result in Section \secb.

\subhead \secza. Preliminaries \endsubhead

Let $W$ be a simply laced Coxeter group with set $S =
\{ s_i : i \in I \}$ of distinguished generators and Coxeter matrix
$( m_{ij})_{i,j \in I}$.  The Coxeter graph, $\Gamma$, of $(W, S)$ has
vertices indexed by $S$, and an edge between distinct vertices 
$s_i$ and $s_j$ if and only if $m_{ij} = 3$.
The basic facts about Coxeter groups needed for
this paper can be found in \cite{{\bf 2}, {\bf 6}}.  We are primarily 
interested in
this paper in the case where $\Gamma$ is a subgraph of a Coxeter graph
of type $E_n$ (for arbitrary $n \geq 6$).

\topcaption{Figure 1} Coxeter graph of type $E_n$
\endcaption
\centerline{
\hbox to 1.500in{
\vbox to 0.500in{\vfill
        \includegraphics{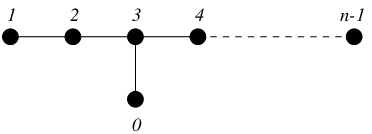}
}
\hfill}
}

The full subgraph of $\Gamma$ omitting vertex $0$ is called the Coxeter
graph of type $A_{n-1}$, and the full subgraph of $\Gamma$ omitting vertex
$1$ is called the Coxeter graph of type $D_{n-1}$.

Denote by $I^*$ the free monoid on $I$. We call the elements of $I$
{\it letters} and those of $I^*$ {\it words}. 
The {\it length} of a word
$\bi \in I^*$ is the number of factors used to express $\bi$
as a product of letters. 
A {\it subword} of a given word $i_1 i_2 \cdots i_n$ (each $i_l \in I$) is
any word of the form $i_p i_{p+1} \cdots i_q$ where $1\leq p \leq q \leq n$.
Let $\phi : I^* \longrightarrow W$ be the surjective
morphism of monoid structures determined by the equalities $\phi(i) = s_i$
for all $i \in I$. We say that a word $\bi \in I^*$ {\it represents} its
image $\phi(\bi ) \in W$.  If the length of a word $\bi \in I^*$ is as small
as possible over all words representing $w = \phi(\bi ) \in W$, then we
call $\bi$ a {\it reduced word}, or a {\it reduced expression} for $w$.
The {\it length} of $w \in W$, denoted by $\ell(w)$, is the
length of any reduced expression for $w$.

Let $V$ be a real vector space with basis $\{ \ga_i : i \in I \}$ in
one-to-one correspondence with $I$. Let $B$ denote the {\it Coxeter form}
on $V$ associated with $(m_{ij})_{i,j \in I}$.  This is the
symmetric bilinear form on $V$ satisfying $B(\ga_i , \ga_j) =
- \cos{\frac{\pi}{m_{ij}}}$ for all $i,j \in I$. Throughout the paper, we
view $V$ as the underlying space of a particular {\it reflection
representation} of $W$, determined by the equalities $s_i \ga_j = \ga_j -
2B( \ga_j , \ga_i ) \ga_i$ for all $i,j \in I$. Note that $B$ is preserved
by $W$ relative to this representation.

Given $u,v \in V$, we say that $u$ is {\it orthogonal} to $v$, and write
$u \perp v$, if $B(u,v) = 0$.

Define $\Phi = \{ w \ga_i : i \in I\}$. This is the {\it root system} of $W$.
The elements of $\Phi$, i.e., the {\it roots}, all have unit length relative
to $B$.  The basis vectors $\ga_i$ are called {\it simple} roots.
Let $\Phi^+$ (respectively, $\Phi^-$) denote the set of all roots expressible as a
linear combination of simple roots with nonnegative (respectively, nonpositive)
coefficients. This is the set of {\it positive} (respectively, {\it
negative}) roots.  We have $\Phi = \Phi^+ \cup \Phi^-$ (disjoint).

Let $w \in W$. Denote by $\Phi(w)$ the set of all $\al \in \Phi^+$ such that
$w \al \in \Phi^-$.  The cardinality of $\Phi(w)$ is $\ell(w)$. Given any
reduced expression $i_1 \cdots i_n$ for $w$, we have
$\Phi(w) = \{ r_1 , \ldots , r_n \}$, where $r_1 = \ga_{i_n}$ and
$r_q = s_{i_n} \cdots s_{i_{n-q+2}}(\ga_{i_{n-q+1}})$ for all $q \in
\{ 2, \ldots , n \}$. The sequence $\br = (r_1 , \ldots , r_n)$ is called
the {\it root sequence} of $i_1 \cdots i_n$, or a root sequence {\it for}
$w$. Note that for any $1 \leq q \leq n$, the initial segment
$(r_1, \ldots , r_q)$ of $\br$ is the root sequence of
$i_{n-q+1} \cdots i_n$.  Observe also that every reduced word is
uniquely determined by its root sequence, so that the map $i_1 \cdots i_n
\mapsto \br$ is a bijection from the set of all reduced expressions for
$w$ to the set of all root sequences for $w$.

\subhead \seczc. Braid moves and contractible triples \endsubhead

For any letters $i, j \in I$ and any positive integer $q$, we define
$(i,j)_q$ to be the length $q$ word $iji\cdots \in I^*$.

Let $\bi , \bj \in I^*$ and let $i,j,k \in I$ with $i \neq j$. Denote the
length of $\bj$ by $n$. We call the
substitution $\bi (i,j)_{m_{ij}} \bj \rightarrow \bi (j,i)_{m_{ji}} \bj$
a {\it braid move}, qualifying it {\it short} or {\it long} according as
$m_{ij}=2$ or $3$, respectively.  Braid moves can be described in terms of
root sequences \cite{{\bf 4}, Proposition 3.1.1}, as follows. Suppose that the word
$\bi ij \bj$ is reduced, and let $\br = (r_q)$ be its root sequence. Then
$m_{ij}=2$ if and only if $r_{n+1} \perp r_{n+2}$, and when these
conditions hold, $\bi ji \bj$ is reduced and its root sequence $\br'$ can
be obtained from $\br$ by swapping $r_{n+1}$ and $r_{n+2}$. We call the
passage from $\br$ to $\br'$ (assuming $r_{n+1} \perp r_{n+2}$) a {\it short
braid move}. Suppose now that the word $\bi ijk \bj$ is reduced, and let
$\br = (r_q)$ be its root sequence. Then $i=k$ if and only if
$r_{n+1} + r_{n+3} = r_{n+2}$, and when these conditions hold,
$\bi jij \bj$ is reduced and its root sequence $\br'$ is obtainable from
$\br$ by swapping $r_{n+1}$ and $r_{n+3}$.  Here, we call the passage from
$\br$ to $\br'$ (assuming $r_{n+1} + r_{n+3} = r_{n+2}$) a
{\it long braid move}.  Our definitions are such that the bijection at the
end of Section \secza\  is compatible with both long and short
braid moves.

Let $w \in W$.  We shall find useful a theorem of Matsumoto \cite{{\bf 8}}
and Tits \cite{{\bf 10}} (stated in \cite{{\bf 1}, Theorem 3.3.1}), which states that 
every reduced expression for $w$ can
be obtained from any other through a sequence of braid moves. A similar
statement holds for root sequences, on account of the aforementioned
bijection.

We say that two words in $I^*$ are {\it short braid equivalent} if one
can be obtained from the other by a sequence of short braid moves.  There
is a corresponding notion for root sequences.

The {\it height} of any $\al \in \Phi$ is the sum of the
coefficients used to express $\al$ as a linear combination of simple
roots. The following definition comes from \cite{{\bf 7}, \S3.1}.

\definition{Definition \seczc.1}
Let $w \in W$.  We call any subset of $\Phi(w)$ of the form
$\{ \al , \be , \al + \be \}$ an {\it inversion triple} of $w$.  If
$T$ is an inversion triple of $w$ such that there is a root sequence
for $w$ in which the elements of $T$ appear consecutively (in some order),
then we say that $T$ is {\it contractible}. We denote by $N(w)$ the
number of contractible inversion triples of $w$, and by
$\widetilde{N}(w)$ the number of roots $\al$ such that $\al$ is the
highest root of at least one contractible inversion triple of $w$.
\enddefinition

For brevity, we usually write ``contractible triple" instead of
``contractible inversion triple".  Observe that $\widetilde{N}(w)
\leq N(w)$ for all $w \in W$.

The following basic properties of inversion triples will be used freely in
the sequel.  (Further details may be found in \cite{{\bf 7}, \S3.1}.)

\proclaim{Lemma \seczc.2}
\item{\rm (i)}
{In the type $A$ setting, every inversion triple is contractible, but
this does not hold in general.}
\item{\rm (ii)}
{Suppose that $( \ldots , \al , \ga , \be , \ldots )$ is a root 
sequence for some $w \in W$, and the following long braid move
can be applied: $( \ldots , \al , \ga , \be , \ldots ) \rightarrow
( \ldots , \be , \ga , \al , \ldots )$.   Then $\{ \al , \be , \ga \}$
is a contractible triple of $w$ with highest root $\ga$.}
\item{\rm (iii)}
{Suppose that $\{ \al , \be , \al + \be \}$ is an inversion
triple of $w$. In every root sequence for $w$, the root $\al + \be$
must appear between $\al$ and $\be$.}
\item{\rm (iv)}
{If $T$ and $T'$ are inversion triples of the same element
of $W$ and $\# (T \cap T') > 1$, then $T=T'$.}
\qed\endproclaim

\definition{Definition \seczc.3 \cite{{\bf 7}, Definition 3.1.6}}
Let $w \in W$.  Suppose that for any pair of intersecting
contractible triples $T$ and $T'$ of $w$, the highest roots in $T$ and
$T'$ agree. Then we say that $w$ is {\it maximally clustered}.  If $W$
has finitely many maximally clustered elements, then we will call $W$
{\it MC-finite}.
\enddefinition

\subhead \seczd. Contracted reduced expressions \endsubhead

We are almost ready to define a key concept for the proof of our main 
result, namely that of a contracted reduced expression for a maximally 
clustered element.

\definition{Definition \seczd.1}
Let $w \in W$ be maximally clustered.  Let $C$ be a collection of
pairwise-intersecting contractible triples of $w$. Suppose that $C$
is nonempty and not properly contained in another set of
pairwise-intersecting contractible triples of $w$. Then we call $C$ a
{\it maximal set of triples} for $w$.
\enddefinition

\definition{Definition \seczd.2}
Let $w \in W$ be maximally clustered, and let $\br$ be a root sequence
for $w$.  Suppose that $C$ is a maximal set of triples for $w$, and
that $$\br = (\ldots , \al_1 , \al_2 , \ldots ,\al_{n-1} , \al_n,
\ga , \be_n , \be_{n-1} , \ldots , \be_2 , \be_1 , \ldots),$$
where $\{ \{ \al_1 , \be_1 , \ga \} , \ldots , \{ \al_n , \be_n , \ga
\} \} = C$. Then we say that $\br$ is {\it contracted for} $C$. We
say that a root sequence for $w$ is {\it contracted} if it is
contracted for every maximal set of triples for $w$.  We call a reduced
expression for $w$ {\it contracted} if its root sequence is contracted.
\enddefinition

\definition{Definition \seczd.3}
Suppose that $\bi \in I^*$ is of the form
$\bi = i_1i_2 \cdots i_ni_{n+1}i_n \cdots i_2i_1$ ($n \geq 1$), where
$i_1, i_2, \ldots , i_{n+1} \in I$ are distinct and where, for each
$1 \leq q \leq n$, there is a unique $q < r \leq n+1$ such that
$m_{i_q i_r} = 3$.  We call $\bi$ a {\it braid cluster}.
\enddefinition

The following result can be used to reduce questions about maximally clustered
elements to questions about braid clusters.

\proclaim{Proposition \seczd.4 (\cite{{\bf 7}, Corollary 4.3.3})}
Let $w \in W$ be maximally clustered.  The following statements
hold:
\item{\rm (i)}{Every reduced expression for $w$ is short braid
equivalent to a contracted reduced expression.}
\item{\rm (ii)}{If $\bi$ is a contracted reduced expression for $w$,
then $\bi = \bi_0 \bc_1 \bi_1 \bc_2 \bi_2 \cdots
\bc_{\widetilde{N}(w)} \bi_{\widetilde{N}(w)}$, where each $\bc_q$ is a
braid cluster, of length $2n_q + 1$ say, and $N(w) = \sum_q n_q$.}
\qed\endproclaim

The next lemma will be useful in the sequel.

\proclaim{Lemma \seczd.5}
\item{\rm (i)}
{Suppose that $w \in W$ is maximally clustered and that $i \in I$
satisfies $\ell(ws_i) < \ell(w)$.  Then $ws_i$ is maximally
clustered.}
\item{\rm (ii)}
{Suppose that $w \in W$ is maximally clustered and that $i \in I$
satisfies $\ell(s_i w) < \ell(w)$.  Then $s_i w$ is maximally
clustered.}
\item{\rm (iii)}
{If $\bi = i_1 \cdots i_n$ is a reduced expression for a maximally clustered
element of $w$, then every subword of $\bi$ is a reduced expression for
some maximally clustered element.}
\endproclaim

\demo{Proof}
Part (ii) may be seen to hold by extending any root sequence for $s_i w$ to
a root sequence for $w$ by appending a single element.  It then follows that
the contractible triples for $s_i w$ are a subset of the contractible triples
for $w$.

Observe that (iii) follows immediate from (i) and (ii) using an inductive
argument, so it remains to prove (i).

Suppose that $w \in W$ is maximally clustered.
If $\br$ is a root sequence for $ws_i$, then $(\al_i, s_i(\br))$
is a root sequence for $w$.  Combining this with the fact that
$s_i|_V$ is linear and preserves the Coxeter form, we find that
every contractible triple $T$ of $ws_i$ gives rise to a
contractible triple $s_i(T)$ of $w$. Thus, $N(ws_i) \leq N(w)$,
and the inequality is strict if $\al_i$ lies in a contractible
triple of $w$. Further, as $s_i|_V$ is injective, the fact that
$w$ is maximally clustered implies that $ws_i$ is maximally
clustered.
\qed\enddemo

\subhead \secze. Main results \endsubhead

In order to understand the context for the main result of this paper, we now
recall the notions of freely braided elements (introduced by the author
and Losonczy in \cite{{\bf 4}}) and fully commutative elements (introduced by
Stembridge in \cite{{\bf 9}}).

\definition{Definition \secze.1}
An element $w \in W$ is said to be {\it freely braided} if both (a) $w$ is
is maximally clustered and (b) the contractible triples of $w$ are disjoint.

An element $w \in W$ is said to be {\it fully commutative} if all reduced
expressions for $w$ are short braid equivalent.
\enddefinition

In the simply laced case, the fully commutative elements agree with the
$iji$-avoiding elements introduced by Fan \cite{{\bf 3}}.

\proclaim{Lemma \secze.2}
Let $w \in W$.
\item{\rm (i)}{If $w$ is fully commutative, then $w$ is freely braided.}
\item{\rm (ii)}{If $w$ is freely braided, then $w$ is maximally clustered.}
\endproclaim

\demo{Proof}
It was shown in \cite{{\bf 5}, Proposition 1.2.2} that $w$ is fully commutative
if and only if $N(w) = 0$, from which it follows vacuously that $w$ is
freely braided.  This proves (i), and (ii) is immediate from the definitions.
\qed\enddemo

Our main result, which will be proved in \S\secb, is as follows.

\proclaim{Theorem \secze.3}
Let $W$ be a simply laced Coxeter group.  Then the following are equivalent:
\item{\rm (i)}{$W$ has finitely many maximally clustered elements;}
\item{\rm (ii)}{$W$ has finitely many freely braided elements;}
\item{\rm (iii)}{$W$ has finitely many fully commutative elements.}
\endproclaim

\remark{Remark \secze.4}
The implications (i) $\Rightarrow$ (ii) and (ii) $\Rightarrow$ (iii)
are immediate from Lemma \secze.2.
\endremark

\subhead \seca. The operators $\pi_j$ on braid clusters \endsubhead

Throughout \S\seca, we assume that $W$ is a Coxeter group of type $E_n$,
whose Coxeter graph is as shown in Figure 1.

\definition{Definition \seca.1}
Let $\bc = i_1 i_2 \cdots i_{n-1} i_n i_{n-1} \cdots i_2 i_1$ be a 
braid cluster.  We define $\Gamma_\bc$ (the subgraph of
$\Gamma$ induced by $\bc$) to be the full subgraph whose vertices are
the letters $i_j$ appearing in a(ny) reduced expression for $\bc$. 
If $1 \leq j < k \leq n$ and we have $m(i_j, i_k) = 3$, then we orient
$\Gamma_\bc$ by adding an arrow pointing from $i_j$ to $i_k$.
\enddefinition

The next result follows from the above definition and \cite{{\bf 7},
Definition 4.1.1}.

\proclaim{Lemma \seca.2}
The orientation of $\Gamma_\bc$ given by Definition \seca.1 assigns a
unique arrow to each edge.  The orientation is characterized by the fact
that the vertex $i_n$ is the unique sink of $\Gamma_\bc$.
\qed\endproclaim

\proclaim{Lemma \seca.3}
If $\bc$ is a braid cluster, then $\Gamma_\bc$ is a Coxeter graph of
type $A$, $D$ or $E$.
\endproclaim

\demo{Proof}
It follows from Definition \seczd.3 that $\Gamma_\bc$ is connected
and has no circuits.
The remaining assertions follow from the fact that $\Gamma_\bc$ is a subgraph
of $\Gamma$, which is of type $E_n$.
\qed\enddemo

\definition{Definition \seca.4}
Let $\bc$ be a braid cluster with middle letter $i_n$.  We say that $\bc$ is 
{\it normalized} if either (i) $\Gamma_\bc$ has no branch point and $i_n$ is
an extremal vertex of $\Gamma_\bc$ or (ii) $i_n$ is the (unique) branch
point of $\Gamma_\bc$.

We define a {\it weak braid cluster} $\overline{\bc}$ to be any reduced 
expression for a group element represented by a braid cluster $\bc$.
If $w \in W$ is maximally clustered, we define a {\it weakly contracted 
reduced expression} for $w$ to be an expression of the form $$
\overline{\bi} = \bi_0 \overline{\bc}_1 \bi_1 \overline{\bc}_2 \bi_2 \cdots
\overline{\bc}_{\widetilde{N}(w)} \bi_{\widetilde{N}(w)}
$$ such that (a) for each $i$, $\overline{\bc}_i$ and $\bc_i$ are reduced
expressions for the same element and (b) $
\bi = \bi_0 \bc_1 \bi_1 \bc_2 \bi_2 \cdots
\bc_{\widetilde{N}(w)} \bi_{\widetilde{N}(w)}
$ is a contracted reduced expression for $w$.
\enddefinition

\definition{Definition \seca.5}
Let $w \in W$ be maximally clustered, and let $\bi$ be a weakly 
contracted reduced expression for $w$ with $$
\bi = \bi_0 \overline{\bc}_1 \bi_1 \overline{\bc}_2 \bi_2 \cdots
\overline{\bc}_{\widetilde{N}(w)} \bi_{\widetilde{N}(w)}
,$$ where each $\overline{\bc}_q$ is a maximal weak
braid cluster, of length $2n_q + 1$ say, and $N(w) = \sum_q n_q$.
If $\bc_j = \overline{\bc}_j = \bc' i_j \bc''$ is a normalized braid cluster 
with middle letter $i_j$ and induced subgraph $\Gamma_j = \Gamma_{\bc_j}$, then
we define $
\pi_j(\bi) = \bi' \bc'_j \bi''
,$ where $
\bi' = \bi_0 \bc_1 \bi_1 \cdots \bi_{j-1}
,$ $
\bi'' = \bi_j \bc_{j+1} \cdots \bi_{\widetilde{N}(w)}
,$ and $$
\bc'_j = \cases
\bc'' & \text{ if $\Gamma_j$ has a branch point and $\bi' i_j$ is not
reduced},\cr
i_j \bc'' & \text{ otherwise.}\cr
\endcases
$$
\enddefinition

\proclaim{Lemma \seca.6}
Let $w \in W$ be maximally clustered, and let $\bc$ be a maximal normalized
braid cluster in
a weakly contracted reduced expression $\bi = \bi' \bc \bi''$ for $w$.  
\item{\rm (i)}{If $i$ is
a letter of $\bc$ that is not a branch point of the induced subgraph 
$\Gamma_\bc$, then $\bi' i$ and $i \bi''$ are reduced.}
\item{\rm (ii)}{If $s$ is a letter and we can parse 
$\bi' = \bv_1 s \bv_2$ and $\bi'' = \bv_3 s \bv_4$ so that $\bv_2$ and $\bv_3$
each consists of letters commuting with $s$, then $s$ is a branch point and
$s$ occurs in $\bc$.}
\endproclaim

\demo{Proof}
We will prove (i) for the case of $\bi' i$; the other case follows similarly.
Let $\bc = i_1 i_2 \cdots i_{n-1} i_n i_{n-1} \cdots i_2 i_1$, and define
$k$ to be the unique integer $1 \leq k \leq n$ such that $i_k = i$.  Since
$i_k$ is not a branch point, there is at most one vertex $i'$ in $\Gamma_\bc$
such that there is an arrow from $i'$ to $i$.
If no such $i'$ exists, then $i \bc$ is not reduced, and (i) follows
by the exchange condition.

If such an $i'$ does exist, then $i' = i_j$ for some $j < k$, and 
we have $$
i \bc = i \bi_1 i' \bi_2 i \bi_3
,$$ where $\bi_1$ and $\bi_2$ consist entirely of generators distinct from $i$
that commute with $i$.  It follows that the letter $i'$ is involved in
two contractible inversion triples: one involving to the occurrences
of $i, i', i$ listed above, and one involving the occurrences of $i_j, i_n,
i_j$ in $\bc$.  Since $i'$ is the highest root in the first triple but not
in the second, we conclude that $i \bc$ is not maximally clustered.  Since $w$
is maximally clustered, all its subwords must be as well, by Lemma
\seczd.5 (iii), and (i) follows by the exchange condition.

For (ii), let us first suppose that $s$ does not occur in $\bc$.
Since $s \bc s$ is reduced, $\bc$
must contain a letter not commuting with $s$.  Suppose that we have
$1 \leq j < k \leq n$ such that $i_j$ and $i_k$ both fail to commute with
$s$.  Since $\Gamma_\bc$ is a tree (by Lemma \seca.3) and $s$ does not
occur in $s$, $s$ completes the unique path in $\Gamma_\bc$ from $i_j$ to
$i_k$ into a circuit in the Coxeter graph of $W$, which is a contradiction.
We conclude that there is a unique $i_j$ not commuting with $s$.
It now follows from Definition \seczd.3 that $s \bc s$ is a braid
cluster, contradicting the maximality of $\bc$.
This proves that $s$ occurs in $\bc$, and an application of (i) shows that
$s$ must be a branch point.
\qed\enddemo

\proclaim{Lemma \seca.7}
Let $w \in W$ be maximally clustered, let $\bi$ be a weakly 
contracted reduced expression for $w$ in which cluster $\bc_j$ is normalized,
and let $\pi_j(\bi) = \bi' \bc'_\bj \bi''$ be as
in Definition \seca.5, so that $\bi = \bi' \bc_\bj \bi''$.  If we have $$
\pi_j(\bi) = \bu_1 s \bu_2 s \bu_3
$$ for some letter $s$, then $\bu_2$ must contain a letter not commuting
with $s$.
\endproclaim

\demo{Proof}
Suppose that $\bu_2$ consists entirely of letters commuting with $s$.

If the right-hand occurrence of $s$ shown comes from $\bi'$, then this
would contradict the fact that $\bu_1$ is reduced.  

If the right-hand occurrence of $s$ comes from $\bi''$, then the left-hand 
occurrence of $s$ must come from $\bi'$, 
otherwise $\bc'_\bj \bi''$ and $\bc_\bj \bi''$ would
not be reduced.  In this case, we must have $\bi' = \bv_1 s \bv_2$ and 
$\bi'' = \bv_3 s \bv_4$, where $\bv_2$ and $\bv_3$ consist of generators
commuting with $s$, which implies that $
\bv_1 \bv_2 (s \bc_j s) \bv_3 \bv_4
$ is a reduced expression for $w$, and that $\bc_j$ contains a generator
not commuting with $s$.  If $\bc_j$ and $\bc'_j$ contain the
same set of generators (disregarding multiplicities), then $\bc'_j$ contains
a generator not commuting with $s$ and the conclusion follows.  The only
other possibility is that we are in the case $\bc'_j = \bc''$ of Definition
\seca.5, and the preceding argument still works unless $s$ is adjacent to
the branch point of $\Gamma_j$.  Since all elements of the Coxeter graph
$\Gamma$ adjacent to the branch point lie in $\Gamma_j$, it follows that
$s$ occurs in the braid cluster $\bc_j$, but that $s$ is not itself a branch
point, and this contradicts Lemma \seca.6 (ii).

We may now assume that the rightmost occurrence of $s$ shown lies in the 
subword $\bc'_\bj$.  Since $\bc'_\bj$ consists of distinct generators, the
leftmost occurrence of $s$ must come from $\bi'$, and $w$ has a reduced
expression of the form $
\bv s \bc_j \bi''
.$  Lemma \seca.6 (i) now forces $s$ to be the branch point of $\Gamma_j$.
Since $\bv s$ is reduced, $\bi' s$ cannot be reduced, and 
Definition \seca.5 shows that we are in the case
$\bc'_j = \bc''$, and this is a contradiction because $s$ does not occur in
$\bc''$.
\qed\enddemo

\proclaim{Lemma \seca.8}
Let $w \in W$ be maximally clustered, and let $\bi$ be as 
in Lemma \seca.7, with $\pi_j(\bi) = \bi' \bc'_\bj \bi''$ 
and $\bi = \bi' \bc_\bj \bi''$.  
Suppose that $$
\pi_j(\bi) = \bu_1 s \bu_2 t \bu_3 s \bu_4
$$ for some noncommuting letters $s$ and $t$, and that $\bu_2$ and $\bu_3$
consist of generators commuting with $s$.  Then the indicated occurrences
of $s$ and $t$ come from the same weak braid cluster $\bc_i$ (where $i \ne j$).
\endproclaim

\demo{Note}
Note that $\bi$ itself satisfies the condition claimed for $\pi_j(\bi)$,
by properties of weakly contracted reduced expressions.
\enddemo

\demo{Proof}
Assume for a contradiction that the indicated occurrences of $s$ and $t$
do not all come from the same braid cluster.

If the right-hand occurrence of $s$ shown comes from $\bi'$, then the 
conclusion follows from the observations that $\bi'$ is maximally clustered 
and its braid clusters are a subset of those of $\bi$.

If the right-hand occurrence of $s$ comes from $\bi''$, then the left-hand 
occurrence of $s$ must come from $\bi'$, as in the proof of Lemma \seca.7.
Assume that the two occurrences of $s$ do not come from the same cluster.
By the Note above, these occurrences must be separated in $\bi$ by at least two
generators not commuting with $s$, and it follows that one of these two
(an occurrence of $u$, say) must have been deleted by $\pi_j$, and that 
$u = t$.  In turn, this means that we have both $\bi' = \bv_1 s \bv_2$ 
and $\bi'' = \bv_3 s \bv_4$ so that $\bv_2$ and $\bv_3$
each consists of letters commuting with $s$.  By Lemma \seca.6 (ii), $s$
occurs in $\bc$, and $s$ is a branch point of $\Gamma_j$.  This means
that the occurrences of $s$ in $\pi_j(\bi)$ are separated by 
occurrences of three different generators not commuting with $s$, a 
contradiction.

We may now assume that the rightmost occurrence of $s$ comes from $\bc'_j$.
Suppose first that $s$ is not the leftmost letter of $\bc'_j$; in 
particular, $s$ is not a branch point of $\Gamma_j$.  
Let $u$ be the unique letter of $\bc_j$ such that there is an arrow from
$s$ to $u$ in $\Gamma_j$.  Since $u$ lies between the two occurrences of
$s$ in $\pi_j(\bi)$, we must have $u = t$, so that in particular, $t$
lies in $\bc'_j$.
This implies that $\bi' = \bv s \bv'$, where $\bv'$ consists
of generators commuting with $s$, which contradicts Lemma \seca.6 (i).

From now on, we may assume that $s$ is the leftmost letter of $\bc'_j$, 
and therefore that $\bi' = \bv_1 s \bv_2 t \bv_3$, where $\bv_2$ and
$\bv_3$ consist of generators commuting with $s$.

Suppose first that we are in the case $\bc'_j = \bc''$ of Definition \seca.5, 
and let $u \ne s$ be the the middle letter of $\bc_j$.  Then $u$ is a branch
node, $s$ is not, and there are precisely two generators in $S$, 
$u$ and $u'$, that
do not commute with $s$.  We have $\bc_j = \bc' sus \bc''$, where $\bc'$
contains no occurrences of $u$, and the expression for $\bi$ in the previous
paragraph shows that $\bi' \bc_j$ is short braid equivalent to a reduced
expression containing the subword $stsus$.  If $t = u$, this is not reduced,
and if $t = u'$, this is not maximally clustered, as the middle $s$ is 
involved in two contractible inversion triples without corresponding to
the highest root of either.  

We must therefore have $\bc'_j = s \bc''$ in Definition \seca.5, so that $s$
is the middle letter of the cluster $\bc_j$.  
We also have $\bi' = \bv_1 t \bv_2$, where
every letter of $\bv_2$ commutes with $t$.

Suppose first that $s$ is a branch point of $\Gamma_j$.  In this case, $t$ 
is not a branch point, but since $t$ is
adjacent to $s$ (both in $\Gamma$ and in $\Gamma_j$), it must be the case
that $t$ occurs in $\bc_j$.  
The fact that the expression $\bi' = \bv_1 t \bv_2$ is reduced
now contradicts Lemma \seca.6 (i).  

The only other possibility is that $s$ is an endpoint of
$\Gamma_j$.  In this case, there is a unique generator $u$ appearing in
$\bc_j$ such that $su \ne us$, so that we have $\bc_j = \bv_3 usu \bv_4$, 
where every letter of $\bv_3$ commutes with $s$.
Since $\bc_j$ is normalized, $u$ is not a branch point of $\Gamma_j$,
and Lemma \seca.6 (i) implies that $\bi' u$ is reduced,
meaning that $t \ne u$ and $t$ does not occur in $\bc_j$.  Since the 
Coxeter graph contains no circuits
and $t$ is adjacent to $s$, Definition \seczd.3 shows that 
$t$ commutes with every generator in $\bc_j$ apart from $s$.  It follows
that $\bi$ is short braid equivalent to a reduced expression with
$stusu$ as a subword, but $stusu = stsus$ is not maximally clustered,
a contradiction.
\qed\enddemo

\proclaim{Lemma \seca.9}
Let $w \in W$ be maximally clustered, and let $\bi$ be as 
in Lemma \seca.7, with $\pi_j(\bi) = \bi' \bc'_\bj \bi''$ 
and $\bi = \bi' \bc_\bj \bi''$.  
Then $\pi_j(\bi)$
is a reduced expression for some element $w' \in W$, and every reduced 
expression for $w'$ is short braid equivalent to $\pi_j(\bk)$ for some
weakly contracted reduced expression $\bk$ for $w$.
\endproclaim

\demo{Proof}
We apply the well known theorem on braid moves due to Matsumoto and Tits,
as stated in \cite{{\bf 1}, Theorem 3.3.1}, starting with the expression 
$\pi_j(\bi)$.

By Lemma \seca.7, it is not possible to apply a sequence of short braid
relations followed by removal of a consecutive pair $ss$ to $\pi_j(\bi)$.
Suppose instead that we apply a sequence of short braid relations followed
by a long braid relation to $\pi_j(\bi)$.  By Lemma \seca.8, this long
braid relation involves three letters from a weak braid cluster $\bc_i$,
where $i \ne j$.  Applying this long braid relation to $\bi$ results in another
weakly contracted reduced expression in which the cluster $\bc_j$ is still
normalized.  It follows that application of the long braid relation commutes
(up to short braid equivalence) with the map $\pi_j$, as required.  The 
procedure may now be iterated.

By \cite{{\bf 1}, Theorem 3.3.1}, any reduced expression for $\bi$ can be obtained
by application of short and long braid relations, together with excision
of any consecutive pairs $ss$.  Iterating the procedure of the last paragraph
shows that we never have an opportunity to remove a pair of the form $ss$,
which means that $\pi_j(\bi)$ is reduced.  The claim about short braid
equivalence also follows immediately from the proof in the previous paragraph.
\qed\enddemo

\proclaim{Lemma \seca.10}
Let $w \in W$ be maximally clustered, and let $\bi$ be a weakly 
contracted reduced expression for $w$.  Then $\pi_1(\bi)$ is a reduced
expression for $w' \in W$, where $w'$ is maximally clustered and
$\widetilde{N}(w') = \widetilde{N}(w) - 1$.
\endproclaim

\demo{Proof}
We know from Lemma \seca.9 that $w'$ is reduced.  If $T$ is a contractible
inversion triple for $w'$, then there is a reduced expression $\bv sts \bv'$
for $w'$ in which $s$ and $t$ are noncommuting generators and the subword
$sts$ corresponds to the triple $T$.  By Lemma \seca.9, there is a weakly
contractible reduced expression $$
\overline{\bi} = \bi_0 \overline{\bc}_1 \bi_1 \overline{\bc}_2 \bi_2 \cdots
\overline{\bc}_{\widetilde{N}(w)} \bi_{\widetilde{N}(w)}
$$ for $w$ such that $\bv sts \bv'$ is short braid equivalent to 
$\pi_1(\overline{\bi})$.  By Lemma \seca.8, the indicated occurrences
of $s$ and $t$ come from the same weak braid cluster, $\overline{\bc}_i$,
where $i > 1$.  Since $w$ and $w'$ both have reduced expressions ending
in $$
\bi_1 \overline{\bc}_2 \bi_2 \cdots
\overline{\bc}_{\widetilde{N}(w)} \bi_{\widetilde{N}(w)}
,$$ it follows that the triple $T$ is also a contractible inversion triple
for $w$, and conversely that every contractible inversion triple of $w$
not arising from the cluster $\bc_1$ is a contractible inversion triple
of $w'$.  The assertions now follow.
\qed\enddemo

\proclaim{Lemma \seca.11}
Let $w \in W$ be maximally clustered, and let $\bi$ be a
contracted reduced expression for $w$.  Applying the
operator $\pi_1$ $\widetilde{N}(w)$ times to $\bi$, we obtain a reduced
expression for a fully commutative element.
\endproclaim

\demo{Proof}
By Lemma \seca.10, the element obtained is a maximally clustered element
$w'$ with $\widetilde{N}(w') = 0$.  Since $w'$ has no contractible inversion
triples, it follows that no reduced expression for $w'$ can have a subword
of the form $iji$, so that $w'$ is fully commutative.
\qed\enddemo

\subhead \secb. Proof of Theorem \secze.3 \endsubhead

A key ingredient of the proof of our main result is the following

\proclaim{Proposition \secb.1}
Let $W$ be a simply laced Coxeter group.  Then the following are equivalent:
\item{\rm (i)}{$W$ has finitely many fully commutative elements;}
\item{\rm (ii)}{the connected components of the Coxeter graph of $W$ are
subgraphs of the Coxeter graph of type $E_n$ (see Figure 1).}
\endproclaim

\demo{Proof}
This is a special case of Stembridge's result \cite{{\bf 9}, Theorem 4.1}.
\qed\enddemo

\proclaim{Theorem \secb.2}
A Coxeter group $W$ of type $E_n$ (for any $n \geq 6$) has finitely many
maximally clustered elements.
\endproclaim

\demo{Proof}
For each maximally clustered element $w \in W$, choose a contracted reduced
expression $\bi(w)$.  By Lemma \seca.11, there is a map $f$ from the set
$\{\bi(w) : w \in W\}$ to the set of reduced expressions for fully commutative
elements of $W$ given by $$
f(\bi(w)) = \pi_1^{\widetilde{N}(w)}(\bi(w))
.$$  By Proposition \secb.1, there are only finitely many fully commutative 
elements of $W$, and since each one has only finitely many reduced
expressions, the problem reduces to showing that the fibres of $f$ are finite.

If we write $$
\bi(w) = \bi_0 \bc_1 \bi_1 \bc_2 \bi_2 \cdots
\bc_{\widetilde{N}(w)} \bi_{\widetilde{N}(w)}
,$$ then $$
f(\bi(w)) = \bi_0 \bc'_1 \bi_1 \bc'_2 \bi_2 \cdots
\bc'_{\widetilde{N}(w)} \bi_{\widetilde{N}(w)}
,$$ where $\bc'_i$ is as in Definition \seca.5.  Given any reduced expression
in the image of $f$, there are only finitely many ways to select the
subwords $\bc'_i$, and for each such subword, there are at most two possible
$\bc_i$ that could have given rise to it.  It follows that the fibres of
$f$ are finite, as required.
\qed\enddemo

\demo{Proof of Theorem \secze.3}
By Remark \secze.4, it only remains to prove the implication (iii)
$\Rightarrow$ (i), and the problem immediately reduces to the case where
the Coxeter graph of $W$ is connected.
Suppose that $W$ has finitely many fully commutative elements.  By Proposition
\secb.1, the Coxeter graph of $W$ is a subgraph of the Coxeter graph of
type $E_n$, so it is enough to deal with the case of $W$ being of 
type $E_n$.  The result now follows from Theorem \secb.2.
\qed\enddemo

\head Acknowledgements \endhead

I thank J. Losonczy for some very helpful correspondence during the 
preparation of this paper.

\leftheadtext{} \rightheadtext{}

\end